\documentclass[12pt]{article} 
\usepackage{amsmath,amsthm, amssymb} 
\usepackage{amssymb,latexsym}

\newtheorem{theorem}{Theorem}

\newtheorem{lemma}{Lemma}

\begin{document}

\baselineskip=17pt 

\title{\bf On the equation $x_1^2 + x_2^2 + x_3^2 + x_4^2 = N$
with variables such that $x_1 x_2 x_3 x_4 + 1$ is an almost-prime}

\author{T. L. Todorova \and D. I. Tolev \footnote{Supported by Sofia University Grant 99/2013}}

\date{}
\maketitle

\begin{abstract}
We consider Lagrange's equation 
$x_1^2 + x_2^2 + x_3^2 + x_4^2 = N$, where $N$ is a sufficiently large and odd integer, 
and prove that it has a solution in natural numbers $x_1, \dots, x_4 $ such that $x_1 x_2 x_3 x_4 + 1$ 
has no more than 48 prime factors.

\medskip


Mathematics Subject Classifcation (2010): 11P05, 11N36.
\end{abstract}

\section{Introduction and statement of the result}

\indent

In 1770 Lagrange proved that for any positive integer $N$
the equation
\begin{equation} \label{10}
 x_1^2 + x_2^2 + x_3^2 + x_4^2 = N
 %\qquad\qquad(10)
\end{equation}
has a solution in integers $x_1, \dots, x_4$ and later 
Jacobi found an exact formula for the number of the solutions
(see \cite[Ch. 20]{Hardy-Wright}).
A lot of researchers studied the equation \eqref{10} for solvability in integers 
satisfying additional conditions. There is a hypothesis stating that if $N$ is sufficiently large 
and $N \equiv 4 \pmod{24}$
then \eqref{10} has a solution in primes. 
This hypothesis has not been proved so far, but several approximations to it have been established.

%\bigskip

Greaves~\cite{Grieves}, Plaksin~\cite{Plaksin}, Shields~\cite{Shields} and Kowalchik~\cite{Kowalchik} 
considered \eqref{10} with two prime and two integer variables. Br\"udern and Fouvry~\cite{Bru-Fouv}, 
Heath-Brown and Tolev~\cite{HB-Tolev}, Tolev~\cite{Tolev}, Yinhchun Cai~\cite{Cai} studied \eqref{10}
with multiplicative restrictions imposed on all of the variables ---
with four almost-primes or with one prime and three almost-primes.
(We say that the integer $n$ is an almost-prime of order $r$ if $n$ has at most $r$ prime factors, counted with the 
multiplicity. We denote by $\mathcal P_r$ the set of all almost-primes of order $r$).

%\bigskip

Yinhchun Cai established in \cite{Cai} the solvability of \eqref{10} in:

%\bigskip

---  $x_1$ prime and $x_2, x_3, x_4 \in \mathcal P_{42}$;

---  $x_1$ prime and $x_2, x_3, x_4 $ satisfying $x_2x_3x_4 \in \mathcal P_{121}$;

---  $x_1, x_2, x_3, x_4 \in \mathcal P_{13}$;

---  $x_1, x_2, x_3, x_4 $ satisfying  $x_1x_2x_3x_4 \in \mathcal P_{41}$.

%\bigskip

We should also mention the result of Blomer and Br\"udern~\cite{Blo-Bru} which states that every
sufficiently large integer, satisfying certain natural congruence conditions, can be represented 
in the form $x_1^2 + x_2^2 + x_3^2$ with integers $x_1, x_2, x_3 \in \mathcal P_{521}$.
Later L\"u Guangshi~\cite{LuGuangshi} considered the same problem, but with integers such that 
$x_1x_2x_3 \in \mathcal P_{551}$. Obviously from these results one obtains information about
the solvability of \eqref{10} in three almost-prime variables and one variable of any nature.

%\bigskip

Having in mind the results mentioned above one may consider
the following problem. 
For a given polynomial $f \in \mathbb Z [x_1, \dots, x_4]$ to study the arithmetical 
properties of the integers $f(x_1, \dots, x_4)$, 
where $x_1, \dots, x_4$ are solutions of \eqref{10} and, in particular, to 
study the solvability of  \eqref{10} in integers 
$x_1, \dots, x_4$ such that $f(x_1, \dots, x_4)$ is an almost-prime of a given order.

%\bigskip

In the present paper we consider the polynomial $f = x_1 x_2 x_3 x_4 + 1$ and prove the following:
\begin{theorem} \label{theorem1}
Suppose that $N$ is a sufficiently large odd integer. Then the equation \eqref{10}
has a solution in natural numbers $x_1, \dots, x_4 $ such that $x_1 x_2 x_3 x_4 + 1$ 
has no more than 48 prime factors. The number of such solutions is greater than 
$ \frac{c N}{\log N}$ for some constant $c>0$. 
\end{theorem}

%\bigskip

A similar result holds if $N$ is even, but having a large odd divisor. 
(If $N$ is a power of $2$ then, according to the Jacobi theorem \cite[Ch. 20]{Hardy-Wright}, the equation
\eqref{10} has exactly $24$ solutions in integers and in this case our method does not work). 
Using the method of the proof one may study 
this problem with an arbitrary polynomial $f \in \mathbb Z[x_1, \dots, x_4]$, satisfying certain natural conditions.

%\bigskip

In the present paper we use the following notations.

%\bigskip

We denote by $N$ a sufficiently large odd integer. 
Letters $a, b, k, l, m, n, v$ are always integers, $q$ is a natural number and
$p$ is always a prime number. 
By $(n_1, \dots, n_k)$ we denote the greatest common divisor of $n_1, \dots, n_k$.
If $q \in \mathbb N$ and $a \in \mathbb Z$, $(a, q)=1$ then we denote by 
$\overline {(a)_q}$ the inverse of $a$ modulo $q$, i.e. the solution of the 
congruence $ax \equiv 1 \pmod{q}$. If the value of the modulus is clear from the 
context then we write $\overline{a}$ for simplicity. 
If $p^l \mid m$, but $p^{l+1} \nmid m$ then we write $p^l \parallel m$.
We denote by $\vec{n}$ four dimensional vectors and let
\begin{equation} \label{290}
  |\vec{n}| = \max (|n_1|, \dots, |n_4|) .
 %\qquad\qquad(290)
\end{equation}
For an odd  $q$ we denote by $\left( \frac{ \cdot}{q} \right)$ the Jacobi symbol.
As usual $\mu(q)$ is the M\"obius function,
$\varphi(q)$ is the Euler function and
$\tau(q)$ is the number of positive divisors of $q$.
Sometimes we write $a \equiv b \; (q)$ as an abbreviation of $a \equiv b \pmod{q}$.
We write $\sum_{x \; (q)}$ for a sum over a complete system of residues modulo $q$
and respectively  $\sum_{x \; (q)^*}$ is a sum over a reduced system of 
residues modulo $q$.
We also denote $e(t)= e^{2 \pi i t}$.

%\bigskip

We use Vinogradov's notation $A \ll B$, which is equivalent to $A = O(B)$. If we have simultaneously 
$A \ll B$ and $B \ll A $ then we write $A \asymp B$.
By $\varepsilon$ we denote arbitrarily small positive number, which is not the same in different
formulas. The constants in the $O$-terms and
$\ll$-symbols are absolute or depend on $\varepsilon$.

\section{Results about exponential and character 
sums and integrals}

\indent

Consider first some classical exponential sums.

%\bigskip

The Gauss sum is defined by
\begin{equation} \label{250}
  G(q, m, n) = \sum_{x \; (q)} e \left( \frac{m x^2 + n x}{q} \right) .
 %\qquad\qquad(250)
\end{equation}
We denote also
\begin{equation} \label{250.5}
  G(q, m) = G(q, m, 0) .
 %\qquad\qquad(250.5)
\end{equation}

%\bigskip

The Gauss sum has the following properties.

%\bigskip

If $(q, m)=d$ then
\begin{equation} \label{251}
  G(q, m, n) = 
  \begin{cases}
     d \,  \displaystyle {G \left( \frac{q}{d}, \frac{m}{d}, \frac{n}{d} \right)  } \qquad & \text{if} \qquad d \mid n , \\
     0 & \text{otherwise} .
  \end{cases}
 %\qquad\qquad(251)
\end{equation}

%\bigskip

For any $q$ we have
\begin{equation} \label{253}
  G(q, 1) = \frac{1 + i^{-q}}{1 + i^{-1}} \sqrt{q} .
 %\qquad\qquad(253)
\end{equation}

%\bigskip

If $(q, 2m)=1$ then 
\begin{equation} \label{252}
  G(q, m, n) = e \left( - \frac{\overline{4m}\,  n^2}{q} \right) \left( \frac{m}{q} \right) G(q, 1) .
 %\qquad\qquad(252)
\end{equation}

%\bigskip

If $2 \nmid m$ and $k \ge 2$ then
\begin{equation} \label{254}
  G(2^k, m, n) 
  = 
  \begin{cases}
  e \left( - \frac{\overline{m} \; (n/2)^2}{2^k} \right) \, 
  2^{\frac{k + 1}{2}} \, c(m, k)  \quad & \text{if} \quad 2 \mid n ,
   \\
   0  & \text{otherwise} .
  \end{cases}
 %\qquad\qquad(254)
\end{equation}
where 
\begin{equation} \label{254.5}
  c(m, k) = 
   \begin{cases} \frac{1 + i^m}{\sqrt{2}} \quad & \text{if}  \quad 2 \mid k  ,
    \\
      e \left( \frac{m}{8} \right) & \text{if}  \quad2 \nmid k.
    \end{cases}
 %\qquad\qquad(254.5)
\end{equation}
In particular, we have 
\begin{equation} \label{255}
   c(m, k)^4 = -1 .
 %\qquad\qquad(255)
\end{equation}

%\bigskip

If $p > 2$ is a prime then for any $m$ we have
\begin{equation} \label{256}
   \sum_{x \; (p)} \left( \frac{x}{p} \right) e \left( \frac{mx}{p} \right)
   = \left( \frac{m}{p} \right) G(p, 1) .
 %\qquad\qquad(256)
\end{equation}
The proofs of formulas \eqref{251} -- \eqref{256} available in \cite[Sec. 6]{Esterman}
and \cite[Ch. 7]{Hua}.

%\bigskip

For $\vec{n} \in \mathbb Z^4$ we denote
\begin{equation} \label{260}
    G(q, m , \vec{n} ) = \prod_{j = 1}^4 G(q, m, n_j) .
 %\qquad\qquad(260)
\end{equation}

%\bigskip

The Kloosterman sum is defined by
\begin{equation} \label{545}
  K(q, m, n) = \sum_{x \; (q)^*} e \left( \frac{m x + n \overline{x}}{q} \right) .
  %\qquad\qquad(545)
\end{equation}
We use A.Weil's bound
\begin{equation} \label{545.5}
  |K(q, m, n)| \le \tau(q) \, q^{\frac{1}{2}} \, (q, m, n)^{\frac{1}{2}} .
  %\qquad\qquad(545.5)
\end{equation}
A proof of \eqref{545.5} is available in \cite[Ch. 11]{IwKo}.

%\bigskip

The Ramanujan sum is defined by
\begin{equation} \label{545.6}
  c_q(m) = K(q, m, 0) 
  %\qquad\qquad(545.6)
\end{equation}
and we have
\begin{equation} \label{545.7}
  c_q(m) = \frac{\mu \left( \frac{q}{d} \right) }{\varphi \left( \frac{q}{d} \right)}  \, \varphi(q)  , \qquad \text{where} \qquad d = (q, m) .
  %\qquad\qquad(545.7)
\end{equation}
For a proof see \cite[Ch. 16]{Hardy-Wright}.

%\bigskip

We need also an estimate for a special character sum. Suppose that $p>2$ is a prime and
$f \in \mathbb F_p [x]$  is a polynomial of degree $k$, which 
is not of the form $c g^2(x)$, where $c$ is a constant and $g \in \mathbb F_p[x]$.
Then we have
\begin{equation} \label{27}
  \left| \sum_{x \; (p)} \left( \frac{f(x)}{p} \right) \right|
  \le (k-1) \sqrt{p} .
  %\qquad\qquad(27)
\end{equation}
For a proof we refer the reader to \cite[Ch. 11]{IwKo}.

%\bigskip

Consider now some exponential integrals.

%\bigskip

We take the infinitely many times differentiable function
\begin{equation} \label{30}
  \omega_0(t) = 
  \begin{cases}
  \exp \frac{1}{ \left(t - \frac{1}{2} \right)^2 - \frac{1}{16} } \quad & \text{for} \quad 
  t \in \left( \frac{1}{4}, \frac{3}{4} \right) , \\
   0 & \text{otherwise} 
  \end{cases}
 %\qquad\qquad(30)
\end{equation}
and define
\begin{equation} \label{91}
      J(\gamma, u ) = \int_{- \infty}^{\infty} \omega_0 (x) \,
      e (\gamma x^2 + u x) \, d x  .
 %\qquad\qquad(91)
\end{equation}
We have
\begin{equation} \label{92}
      J(\gamma, u ) \ll \min (1, |\gamma|^{- \frac{1}{2}}) .
 %\qquad\qquad(92)
\end{equation}
A proof can be found for example in \cite[Ch. 1]{Karatsuba}.

%\bigskip

For $\vec{u} \in \mathbb R^4$ we define
\begin{equation} \label{280}
  J (\gamma, \vec{u}) = \prod_{j=1}^4 J(\gamma, u_j) 
 %\qquad\qquad(280)
\end{equation}

%\bigskip

We specify the constant $\varkappa$ by
\begin{equation} \label{90}
    \varkappa = \int_{- \infty}^{\infty} e (- \gamma ) \,
      J(\gamma, \vec{0} )  \, d \gamma .
 %\qquad\qquad(90)
\end{equation}
Using the standard technique of the circle method (see for example \cite[Ch. 11]{Karatsuba}) 
one can establish that
\begin{equation} \label{93}
      \varkappa > 0 .
 %\qquad\qquad(93)
\end{equation}

%\bigskip

If $\vec{u} \in \mathbb R^4$ and $|\vec{u}| > 0$
(see \eqref{290} for the definition of $|\vec{u}|$) then we have
\begin{equation} \label{93.5}
      \int_{- \infty}^{\infty} |J(\gamma, \vec{u})| \, d \gamma \ll |\vec{u}|^{-1 + \varepsilon} .
 %\qquad\qquad(93.5)
\end{equation}
The proof of this estimate is available in \cite[Lemma 10]{HB-Tolev}.

\section{Proof of the theorem}

\subsection{Beginning of the proof}

\indent

We denote
\begin{equation} \label{20}
  P = \sqrt{N}  
 %\qquad\qquad(20)
\end{equation}
and let
\begin{equation} \label{40}
  \omega(t) = \omega_0 \left( \frac{t}{P} \right) ,
 %\qquad\qquad(40)
\end{equation}
where $\omega_0(t)$ is defined by \eqref{30}.

%\bigskip

Suppose that $\eta > 0$ is a constant, which will be specified later and let
\begin{equation} \label{21}
\qquad z = N^{\eta} ,  \qquad P(z) = \prod_{2 < p < z} p .
 %\qquad\qquad(21)
\end{equation}
Consider the sum
\begin{equation} \label{50}
  \Gamma = \sum_{\substack{
  x_1^2 + \dots + x_4^2 = N \\ ( x_1 x_2 x_3 x_4 + 1 , P(z) ) = 1 }} 
  \omega(x_1) \dots \omega (x_4) .
 %\qquad\qquad(50)
\end{equation}
If we prove the inequality
\begin{equation} \label{60}
  \Gamma \gg \frac{N}{\log N} 
 %\qquad\qquad(60)
\end{equation}
then we will establish that there is a constant $c>0$ such that
the equation \eqref{10} has at least 
$\frac{c N}{\log N}$ solutions satisfying $(x_1x_2x_3x_4 + 1, P(z))=1$
and such that $x_1 x_2 x_3 x_4 + 1 \asymp N^2$.
We also note that $2 \nmid x_1 x_2 x_3 x_4 + 1$ because in the opposite case
we would have $2 \nmid x_j$ for all $j$ which would imply $2 \mid N$, but this contradicts our assumption.
Hence for every such solution the integer $x_1x_2x_3x_4 + 1$ 
does not have prime factors less than $z$ and therefore
this integer has at most $2/{\eta}$ prime factors.
So, to prove Theorem~\ref{theorem1}, we 
have to choose
$\eta = \frac{1}{24} - \omega $, where $\omega > 0$ is a sufficiently small constant,
and to establish \eqref{60}.

%\bigskip

To find the lower bound \eqref{60} we apply the linear sieve
and that is why we need information about the sums
\begin{equation} \label{70}
  F(N, d) = 
   \sum_{\substack{
  x_1^2 + \dots + x_4^2 = N \\ 
   x_1 x_2 x_3 x_4 + 1 \equiv 0 \; (d) }} 
  \omega(x_1) \dots \omega (x_4) 
 %\qquad\qquad(70)
\end{equation}
where $d$ is squarefree and odd.
Applying the Kloosterman form of the
Hardy--Littlewood circle method, we find
that for small $d$ the sum \eqref{70} can be approximated  by the quantity
\begin{equation} \label{85}
   M(N, d) = \varkappa \, N \, a(N) \, \Psi(N, d) ,
 %\qquad\qquad(85)
\end{equation}
where the terms in the right-hand side of \eqref{85}
are defined as follows.

%\bigskip

The constant $\varkappa$ is given by \eqref{90}.

%\bigskip

Further
\begin{equation} \label{100}
    a(N) = \prod_{p>2} \left( 1 + \frac{1}{p} \right) \left( 1 - \frac{1}{p^{1 + \xi_p(N)}} \right) ,
 %\qquad\qquad(100)
\end{equation}
where $\xi_p (N)$ is the non-negative integer defined by
\begin{equation} \label{110}
    p^{\xi_p(N)} \parallel N .
 %\qquad\qquad(110)
\end{equation}

%\bigskip

Next we have
\begin{equation} \label{120}
    \Psi(N, d) = 
    \frac{\alpha(N, d) \, \mathcal L (N, d)}{d^3} ,
 %\qquad\qquad(120)
\end{equation}
where
\begin{equation} \label{130}
    \alpha(N, d) = 
    \prod_{p \mid d} \left( 1 + \frac{1}{p} \right)^{-1} \left( 1 - \frac{1}{p^{1 + \xi_p(N)}} \right)^{-1} 
 %\qquad\qquad(130)
\end{equation}
and where $\mathcal L (N, d) $ is the number of solutions of the system of congruences
\begin{equation} \label{140}
      b_1^2 + \dots + b_4^2 \equiv N \pmod{d} , \qquad b_1 b_2 b_3 b_4 + 1 \equiv 0 \pmod{d} .
 %\qquad\qquad(140)
\end{equation}

%\bigskip

We denote by $R(N, d)$ the error which arises when we approximate $F(N, d)$ by $M(N, d)$, 
that is 
\begin{equation} \label{80}
  F(N, d) = M(N, d) + R(N, d) .
 %\qquad\qquad(80)
\end{equation}

%\bigskip

To prove our theorem we have to study the arithmetic properties of the main term and to estimate 
the error term.

\subsection{An estimate for a special exponential sum}

\indent

An important role in our analysis  plays the exponential sum
\begin{equation} \label{390}
   V_q =  V_q(N, d, v, \vec{b}, \vec{n})  = 
     \sum_{ a \; (q)^* }
  e \left( \frac{a(b_1^2 + \dots + b_4^2 -N) + \overline{a} v}{q} \right)
   G(q, ad^2, 2ad \vec{b} + \vec{n}) ,
     %\qquad\qquad(390)
\end{equation}
where $\vec{b} = \langle b_1, \dots, b_4 \rangle \in \mathbb Z^4$.
It is analogous to the sum, considered in \cite[Sec. 1]{Bru-Fouv}.

%\bigskip

To estimate $ V_q(N, d, v, \vec{b}, \vec{n})$ 
we use the properties of the Gauss sum
and the Kloosterman sum and prove following:
\begin{lemma} \label{lemma1}
Suppose that $N, d, q \in \mathbb N$, $v \in \mathbb Z$, $2 \nmid Nd$, 
$\mu^2(d)=1$ and $\vec{n} = \langle n_1, \dots, n_4 \rangle  \in \mathbb Z^4$,
$\vec{b} = \langle b_1, \dots, b_4 \rangle \in \mathbb Z^4$.
Then we have
\begin{equation} \label{420}
      V_q(N, d, v, \vec{b}, \vec{n}) \ll
      \tau(q) \, q^{\frac{5}{2}}  \, (q, N)^{\frac{1}{2}}  \, (q, N-b_1^2 - \dots - b_4^2)^{\frac{1}{2}}
       \,  (q, d^2)^2 ,
     %\qquad\qquad(420)
\end{equation}
where the constant in the $\ll$-symbol is absolute.
Further, if some of the conditions
\begin{equation} \label{430}
      (q, d) \mid n_j , \qquad j=1, \dots, 4 
     %\qquad\qquad(430)
\end{equation}
do not hold, then $V_q(N, d, v, \vec{b}, \vec{n}) = 0 $.
\end{lemma}

\paragraph{Proof.}
Suppose that $(q, d) \nmid n_j$ for some $j$. It follows from \eqref{251}
that 
\[ 
   G(q, ad^2, 2adb_j + n_j) = 0
\] 
and having in mind \eqref{260} and \eqref{390}
we see that $V_q = 0$.

%\bigskip

From this point onwards we assume that \eqref{430} holds and 
we begin the proof of \eqref{420}.

%\bigskip

First we note that the sum $V_{q}$ 
is multiplicative with respect to $q$ in the following sense:
If $(q', q'')=1$ then we have
\begin{equation} \label{540}
  V_{q' q''}(N, d, v, \vec{b}, \vec{n}) = V_{q'}\left( N, \, q'' d,  \, \overline{(q'' )_{q'}}^{\, 2} v,  \, \vec{b},  \, \vec{n} \right) 
  \,
  V_{q''}\left( N,  \, q' d,  \, \overline{(q' )_{q''}}^{\, 2} v,  \, \vec{b},  \, \vec{n} \right) .
  %\qquad\qquad(540)
\end{equation}
We leave the routine calculations to the reader.

%\bigskip

Having in mind the identity \eqref{540} we see that 
it is enough to estimate
$  V_{p^s}(N, Ad, Bv, \vec{b}, \vec{n}) $ where $A, B \in \mathbb Z$ and $ p \nmid A $.

%\bigskip

Consider first the case $p > 2$, $p \nmid d$. Applying \eqref{252} 
we find
\[
  G \left( p^s, a  A^2 d^2, 2aAd b_j + n_j \right) = 
   e \left( - \frac{ \overline{4aA^2 d^2} n_j^2  + \overline{Ad} n_j b_j + a b_j^2 }{p^s} \right)  
   \; \left( \frac{a}{p^s} \right) \; G (p^s, 1) .
\]
Therefore, using \eqref{253}, \eqref{260}, \eqref{545} and \eqref{390} we get
\begin{equation} \label{550}
  V_{p^s}(N, Ad, Bv, \vec{b}, \vec{n}) = 
  p^{2 s} \; e \left( - \frac{\overline{Ad} (n_1 b_1 + \dots + n_4 b_4) }{p^s} \right) \; 
  K \left( p^s, -N, M \right)
  %\qquad\qquad(550)
\end{equation}
where
\[
   M \equiv Bv - \overline{4A^2 d^2} (n_1^2 + \dots + n_4^2) \pmod{p^s} .
\]
Using \eqref{545.5} and \eqref{550} we find
\begin{equation} \label{552}
   | V_{p^s}(N, Ad, Bv, \vec{b}, \vec{n})|
    \le (s + 1) \, \left( p^s \right)^{\frac{5}{2}} \, (p^s, N)^{\frac{1}{2}} 
    \qquad \text{for} \qquad p \nmid 2d .
  %\qquad\qquad(552)
\end{equation}

%\bigskip

Consider now the case $p \mid d$. Since $d$ is squarefree we may write
\begin{equation} \label{555}
  d = pd', \qquad \text{where} \qquad  p \nmid d' .
  %\qquad\qquad(555)
\end{equation}

%\bigskip

If $s = 1$ then, using \eqref{430}, we see that $ G \left( p, a  A^2 d^2, 2aAd b_j + n_j \right) = p$
and, having in mind \eqref{260}, we find $ G ( p, a  A^2 d^2, 2aAd \vec{b} + \vec{n} ) = p^4$.
Therefore, from \eqref{545} and \eqref{390} it follows that
\[
   V_{p} (N, Ad, Bv, \vec{b}, \vec{n}) = p^4 K(p, b_1^2 + \dots + b_4^2 - N, Bv ) .
\]
Noting that $(p, d^2) =  p $ and using \eqref{545.5} we find
\begin{equation} \label{560}
  |V_{p}(N, Ad, Bv, \vec{b}, \vec{n})| \le 2 p^{\frac{5}{2}} (p, N - b_1^2 - \dots - b_4^2)^{\frac{1}{2}}
   (p, d^2)^2 \qquad \text{for} \qquad p \mid d .
  %\qquad\qquad(560)
\end{equation}

%\bigskip

In the case $s \ge 2$ we use the following observation.
From \eqref{251} it follows that
\[
   G \left( p^s, a  A^2 d^2, 2aAd b_j + n_j \right) = 0 
\]
unless 
\[
  (p^s, aA^2 d^2) \mid 2aAdb_j + n_j .
\]
Since $p \nmid aA^2$ we see that the later condition is equivalent to
\begin{equation} \label{570}
  (p^s, d^2) \mid 2aAdb_j + n_j .
  %\qquad\qquad(570)
\end{equation}
(The last formula implies, in particular, that 
$ p \mid n_j $, but we already know this because of the assumption \eqref{430}).
Hence we may write 
\begin{equation} \label{580}
  n_j = p n_j', \qquad n_j' \in \mathbb Z , \qquad 1 \le j \le 4
  %\qquad\qquad(580)
\end{equation}
because otherwise $V_q = 0$.

%\bigskip

Suppose that $s = 2$. Then from \eqref{570} it follows that
$p^2 \mid  2aAdb_j + n_j $, hence using \eqref{260} we get
\[ 
   G(p^2, a A^2 d^2, 2aAd \vec{b} + \vec{n} ) = p^8 .
\]
Now we take into account \eqref{545} and \eqref{390} to find
\[
   V_{p^2} (N, Ad, Bv, \vec{b}, \vec{n}) = p^8 K(p^2, b_1^2 + \dots + b_4^2 - N, Bv) .
\]
Noting that $(p^2, d^2)= p^2$ and using \eqref{545.5} we find
\begin{equation} \label{590}
  |V_{p^2}(N, Ad, Bv, \vec{b}, \vec{n})| \le 3 (p^2)^{\frac{5}{2}} (p^2, N - b_1^2 - \dots - b_4^2)^{\frac{1}{2}}
   (p^2, d^2)^2 \qquad \text{for} \qquad p \mid d .
  %\qquad\qquad(590)
\end{equation}

%\bigskip

Consider now the case $s \ge 3$.
Having in mind \eqref{555}, \eqref{570} and \eqref{580} we denote
\begin{equation} \label{600}
   \frac{2aA d' b_j  + n_j'}{p} = h_j \in \mathbb Z .
  %\qquad\qquad(600)
\end{equation}
Using that $(p^s, d^2) = p^2$ and applying \eqref{251}
we find
\begin{align}
   G (p^s, a A^2 d^2, 2aAd b_j + n_j) 
    & = 
      p^2 \,
   G (p^{s - 2}, a A^2 d'^2, h_j) 
   \notag \\
   & \notag \\
    & = 
   p^2 \, e \left( - \frac{ \overline{ (4aA^2 d'^2)_{p^{s-2}} } \, h_j^2} {p^{s - 2}} \right) \,
   \left( \frac{a}{p^{s-2}} \right) \, G(p^{s - 2}, 1) .
   \notag
\end{align}
It is obvious that $\overline {(M)_{p^{s - 2}}} \equiv \overline {(M)_{p^s} } \pmod{p^{s -2}} $
for any integer $M$ with $p \nmid M$.
Hence, using the definition of $h_j$ given by \eqref{600}, we find
\[
   G (p^s, a A^2 d^2, 2aAd b_j + n_j) = 
   p^2 e \left( - \frac{ \overline{4a A^2 d'^2} \, n_j'^2 
      + \overline{Ad'} \, n_j' b_j + a b_j^2  }{p^s} \right)
       \left( \frac{a}{p^{s-2}} \right) G(p^{s - 2}, 1) 
\]
where the inverses are already taken modulo $p^s$.
Therefore, using \eqref{253} and \eqref{260} we find
\[
    G (p^s, a A^2 d^2, 2aAd \vec{b} + \vec{n}) 
    = p^{2 s + 4} 
    e \left( - \frac{ \overline{4a A^2 d'^2} \, \sum_{j=1}^4 n_j'^2 
      \, + \, \overline{Ad'} \, \sum_{j=1}^4  n_j' b_j \, + \, a \sum_{j=1}^4  b_j^2  }{p^s} \right) .
\]
From this formula, \eqref{545} and \eqref{390} we find
\[
  V_{p^s} (N, Ad, Bv, \vec{b}, \vec{n})
    = p^{2 s + 4} 
      e \left( - \frac{  \overline{Ad'} \, \sum_{j=1}^4  n_j' b_j  }{p^s} \right) 
      K (p^s, - N, M)
\]
where $M \equiv Bv - \overline{4 A^2 d'^2} \sum_{j=1}^4 n_j'^2 \pmod{p^s}$.
Therefore, using \eqref{545.5} we obtain
\begin{equation} \label{610}
   | V_{p^s} (N, Ad, Bv, \vec{b}, \vec{n})| 
   \le (s + 1) ( p^s )^{\frac{5}{2}} (p^s, N)^{\frac{1}{2}} (p^s, d^2)^2 
   \qquad \text{for} \qquad p \mid d , \quad s \ge 3 .
  %\qquad\qquad(610)
\end{equation}

%\bigskip

Combining \eqref{552}, \eqref{560}, \eqref{590} and \eqref{610} we see that 
for any prime $p \nmid 2A$ and for any positive integer $s$ we have
\begin{equation} \label{620}
   | V_{p^s} (N, Ad, Bv, \vec{b}, \vec{n})| 
   \le \tau(p^s) \, (p^s)^{\frac{5}{2}} \, 
   (p^s, N - b_1^2 - \dots - b_4^2)^{\frac{1}{2}} \, 
   (p^s, N)^{\frac{1}{2}} (p^s, d^2)^2 .
  %\qquad\qquad(620)
\end{equation}

%\bigskip

Consider now the case $p=2$. We have $2 \nmid dN$ and suppose also that $2 \nmid A$.
We shall prove that for all positive integers $s$ we have
\begin{equation} \label{650}
   | V_{2^s} (N, Ad, Bv, \vec{b}, \vec{n})| 
   \le 4 (s + 1) \, (2^s)^{\frac{5}{2}}  .
  %\qquad\qquad(650)
\end{equation}

%\bigskip

From \eqref{390} it follows that \eqref{650} is obvious for $s = 1$.
Suppose now that $s \ge 2$.
From \eqref{254}, \eqref{260} and \eqref{390} we see that $V_{2^s}$ vanishes if
$2 \nmid n_j$ for some $j$. Hence we may assume that
$2 \mid n_j$ for all $j$, so we may write
\[
  n_j = 2 n_j',  \qquad \text{where} \qquad n_j' \in \mathbb Z, \qquad1 \le j \le 4 .
\]
Using these formulas, as well as \eqref{254} and \eqref{255}, it is easy to verify that
\begin{align}
  G(2^s, aA^2 d^2, 2 a A d \vec{b} + \vec{n}) 
     = -  2^{2 s + 2} \; 
   &
  e \left( 
    - \frac{ 
          \overline{Ad} \, (n_1 b_1 + \dots + n_4 b_4) 
          } {2^s} 
              \right) 
              \notag \\
              & \notag \\
     \times   \,       
     &
e \left( 
    - \frac{ \overline{a A^2 d^2} \, ({n_1'}^2 + \dots + {n_4'}^2)} 
           {2^s} 
       \right) 
        \,
         e \left(     
          - \frac{ a (b_1^2 + \dots + b_4^2)
                 }
           {2^s} 
           \right)          .
  \notag
\end{align}
Hence using \eqref{545} and \eqref{390} we find
\begin{align}
    V_{2^s} (N, Ad, Bv, \vec{b}, \vec{n})
    & = 
    -  2^{2 s + 2} \;
    e \left( 
    - \frac{ 
          \overline{Ad} \, (n_1 b_1 + \dots + n_4 b_4) 
          } {2^s} 
          \right) 
          \notag \\
          & \notag \\
    & \qquad \qquad \times      
          K(2^s, -N, Bv - \overline{A^2 d^2} ({n_1'}^2 + \dots + {n_4'}^2) ) .
          \notag
\end{align}
It remains to apply \eqref{545.5} and we prove \eqref{650}.

%\bigskip

From \eqref{540}, \eqref{620} and \eqref{650} we obtain \eqref{420}
and the proof of the lemma is complete.

\hfill
$\square$

\subsection{The error term}

\indent

In this section we study the quantity $R(N, d)$ defined by \eqref{80}.
Recall that $\mathcal L(N, d)$ is the number of $\vec{b} = \langle b_1, b_2, b_3, b_4 \rangle  \in \mathbb Z^4$
satisfying
\begin{equation} \label{182}
      1 \le b_1, b_2, b_3, b_4 \le d , \quad 
      b_1 b_2 b_3 b_4 + 1 \equiv 0 \pmod{d} , \quad
      b_1^2 + b_2^2 + b_3^2 + b_4^2 \equiv N \pmod{d} .
  %\qquad\qquad(182)
\end{equation}
We prove the following
\begin{lemma} \label{theorem1.5}
Suppose that $N$ is sufficiently large, $d$ is squarefree, $2 \nmid dN$ and
\begin{equation} \label{148} 
      d \le N^{\frac{1}{12}} .
 %\qquad\qquad(148)
\end{equation}
Then we have
\begin{equation} \label{149}
      R(N, d) \ll \mathcal L(N, d) \, N^{\frac{3}{4} + \varepsilon} .
 %\qquad\qquad(149)
\end{equation}
\end{lemma}

\paragraph{Proof.}

We write the sum $F(N, d)$ specified by \eqref{70} in the form
\begin{equation} \label{170}
      F(N, d) = \sum_{ \vec{b}  \; : \; \eqref{182} } \Phi ( N, d, \vec{b} ) ,
 %\qquad\qquad(170)
\end{equation}
where the summation in \eqref{170} is taken over $\vec{b} $ satisfying \eqref{182}
and where
\begin{equation} \label{180}
       \Phi (N, d, \vec{b})  = 
       \sum_{\substack{x_1^2 + \dots + x_4^2 = N \\
       x_i \equiv b_i \; (d) , \quad i = 1, \dots, 4}} 
       \omega(x_1) \dots \omega(x_4) .
 %\qquad\qquad(180)
\end{equation}

%\bigskip

We express the sum \eqref{180} as
\begin{equation} \label{190}
       \Phi (N, d, \vec{b})  = 
       \int_{\mathcal I} e(- N\alpha) \prod_{j=1}^4 S_{d, b_j} (\alpha) 
        \, d \alpha ,
 %\qquad\qquad(190)
\end{equation}
where the integration is taken over the interval
%\[
 %$ \mathcal I = \left( \frac{1}{1 + [P]}, 1 + \frac{1}{1 + [P]} \right] $
 $ \mathcal I = \big( (1 + [P])^{-1}, 1 + (1 + [P])^{-1} \big] $
%\]
and
\begin{equation} \label{200}
       S_{d, b} (\alpha) = \sum_{\substack{x \in \mathbb Z \\ x \equiv b \; (d)}} 
       \omega(x) e(\alpha x^2) .       
 %\qquad\qquad(200)
\end{equation}
Using the properties of the Farey fractions (see \cite[Ch. 3]{Hardy-Wright}) 
we represent $\mathcal I$ as an union of disjoint intervals 
in the following way:
\begin{equation} \label{210}
  \mathcal I = \bigcup_{q \le P} \bigcup_{\substack{ 1 \le a \le q \\ (a, q)=1 }}
  \mathcal L(q, a) ,
 %\qquad\qquad(210)
\end{equation}
where 
\[
 \mathcal L(q, a) = \left( \frac{a}{q} - \frac{1}{q(q+q')}, \frac{a}{q} + \frac{1}{q(q+q'')} \right]
\]
and where the integers $q', q''$ are specified by
\begin{equation} \label{220}
  P < q+q', q+q'' \le q + P , \qquad aq' \equiv 1 \pmod{q} , \qquad aq'' \equiv -1 \pmod{q} .
 %\qquad\qquad(220)
\end{equation}

%\bigskip

We apply \eqref{190},  \eqref{210} and change the variable
of integration to get
\begin{equation} \label{230}
  \Phi (N, d, \vec{b}) = \sum_{q \le P} \sum_{\substack{ a=1 \\ (a, q)=1}}^q
   e\left( - \frac{aN}{q}  \right)
    \int_{\mathfrak M(q, a)}
      e\left( - \beta N  \right)  \prod_{j=1}^4 S_{d, b_j} \left( \frac{a}{q} + \beta \right) 
        \, d \beta ,
 %\qquad\qquad(230)
\end{equation}
where
\begin{equation} \label{240}
  \mathfrak M(q, a) = \left(  - \frac{1}{q(q+q')},  \frac{1}{q(q+q'')} \right] .
 %\qquad\qquad(240)
\end{equation}

%\bigskip

Working as in the proof of \cite[Lemma 12]{HB-Tolev}, we find
that for $\beta \in \mathfrak M(q, a)$ we have
\[
  S_{d, b_j} \left( \frac{a}{q} + \beta \right) =
  \frac{P}{dq} e \left( \frac{a b_j^2}{q} \right)
  \sum_{|n| \le d P^{\varepsilon}} e \left( \frac{n b_j}{dq} \right)
    J \left(\beta N, - \frac{nP}{dq} \right) 
    G(q, ad^2 , 2 a b_j d + n)
    + O \left( P^{-A} \right) ,
\]
where $G(q, m, n)$ and $J(\gamma, u)$ are defined respectively by \eqref{250} and \eqref{91},
$A$ is an arbitrarily large constant, $\varepsilon > 0$ is arbitrarily small 
and the constant in the $O$-term depends 
only on $A$ and $\varepsilon$.
We leave the verification of the last formula to the reader.

%\bigskip

Therefore we may write the integrand in \eqref{230} in the form
\begin{align}
   & 
  e (- \beta N) \frac{P^4}{d^4 q^4} e \left( \frac{a(b_1^2 + \dots + b_4^2)}{q} \right) 
   \notag \\
   & \notag \\
   & \qquad \qquad  \times
  \sum_{|\vec{n}| \le d P^{\varepsilon}}
  e\left( \frac{n_1 b_1 + \dots + n_4 b_4}{dq} \right) 
  J \left( \beta N, - \frac{P}{dq} \vec{n} \right)
  G(q, a d^2 , 2ad \vec{b} + \vec{n}) 
   + O \left( P^{-A} \right)  ,
  \notag
\end{align}
where $G(q, m, \vec{n})$ 
and $J(\beta, \vec{u})$ are defined respectively by \eqref{260} and \eqref{280}
and where the meaning of $|\vec{n}|$ is explained in \eqref{290}.

%\bigskip

We substitute the above expression for the integrand in \eqref{230}, 
change the variable $\beta N = \gamma$ and use \eqref{240} to find
\begin{equation} \label{295}
  \Phi(N, d, \vec{b}) = \tilde\Phi(N, d, \vec{b})  +  O(1) ,
 %\qquad\qquad(295)
\end{equation}
where
\begin{align} 
  \tilde\Phi(N, d, \vec{b}) 
   & = 
  \frac{P^2}{d^4}
  \sum_{q \le P} q^{-4} \sum_{\substack{ a=1 \\ (a, q)=1 }}^q
  e \left( \frac{a(b_1^2 + \dots + b_4^2 -N)}{q} \right)
    \sum_{|\vec{n}| \le d P^{\varepsilon}} 
      e \left( \frac{n_1 b_1 + \dots + n_4 b_4}{dq} \right)
    \notag \\
    & \notag \\
    & \qquad \qquad \times
  G(q, ad^2, 2ad \vec{b} + \vec{n}) \int_{\mathfrak N(q, a)} e(- \gamma) 
  J \left( \gamma, - \frac{P}{dq} \vec{n}\right) \, d \gamma  
 %\qquad\qquad(300)
  \label{300}
\end{align}
and where
\begin{equation} \label{320}
  \mathfrak N(q, a) = \left( - \frac{N}{q(q+q')} , \frac{N}{q(q+q'')} \right] .
 %\qquad\qquad(320)
\end{equation}
We note that from \eqref{220} and \eqref{320} follows
\begin{equation} \label{330}
    \left( - \frac{P}{2q}, \frac{P}{2q}  \right] \subset \mathfrak N(q, a) 
    \subset \left[ - \frac{P}{q}, \frac{P}{q}  \right] .
 %\qquad\qquad(330)
\end{equation}
Therefore we may represent the expression in \eqref{300} as
\begin{equation} \label{340}
    \tilde\Phi(N, d, \vec{b}) = \Phi'(N, d, \vec{b}) + \Phi''(N, d, \vec{b}) ,
 %\qquad\qquad(340)
\end{equation}
where in $\Phi'(N, d, \vec{b})$ the integration is taken over 
$  \gamma \in  \left[ - \frac{P}{2q}, \frac{P}{2q}  \right]$
and, respectively, in $\Phi''(N, d, \vec{b})$ we integrate over 
$ \gamma \in  \mathfrak N(q, a)  \setminus \left[ - \frac{P}{2q}, \frac{P}{2q}  \right]$.

%\bigskip

Consider first $\Phi''(N, d, \vec{b})$. We change the order of summation over $a$ and integration over 
$\gamma$. Using \eqref{330} we conclude that in the new expression for $\Phi''$ the domain of
integration is $ \frac{P}{2q} \le |\gamma| \le \frac{P}{q}$ and 
in the domain of summation over $a$ is imposed the additional condition $\mathfrak N(q, a) \ni \gamma$.
The later condition may be expressed using the idea of Kloosterman \cite{Klos}
and an explanation of this method is available also in  \cite[Sec. 3]{HB-1}.

%\bigskip

There exists a function $\sigma (v, q, \gamma)$, defined for $q \le P$, 
$|\gamma| \le \frac{P}{q}$, $- \frac{q}{2} < v \le \frac{q}{2} $, 
integrable with respect to $\gamma$, satisfying
\begin{equation} \label{350}
    |\sigma (v, q, \gamma)| \le (1 + |v|)^{-1}
 %\qquad\qquad(350)
\end{equation}
and also
\begin{equation} \label{360}
    \sum_{- \frac{q}{2} < v \le \frac{q}{2}} 
     e \left( \frac{\overline{a} v}{q} \right)
    \sigma (v, q, \gamma) =
    \begin{cases}
    1  \quad &  \text{if} \quad \gamma \in \mathfrak N(q, a),
     \\
    0  & \text{otherwise} .
    \end{cases}
 %\qquad\qquad(360)
\end{equation}
Hence using  \eqref{330} and \eqref{360}
we may write $\Phi''$ in the form
\begin{align} 
  & \Phi''(N, d, \vec{b}) 
    = 
  \frac{P^2}{d^4}
  \sum_{q \le P} q^{-4}
      \sum_{|\vec{n}| \le d P^{\varepsilon}} 
      e \left( \frac{n_1 b_1 + \dots + n_4 b_4}{dq} \right)
        \int_{  \frac{P}{2q} \le |\gamma| \le \frac{P}{q}} 
         e(- \gamma) 
  J \left( \gamma, - \frac{P}{dq} \vec{n}\right)
         \notag \\
    & \notag \\
    & \qquad \times
      \sum_{\substack{ a=1 \\ (a, q)=1 }}^q
  e \left( \frac{a(b_1^2 + \dots + b_4^2 -N)}{q} \right)
  G(q, ad^2, 2ad \vec{b} + \vec{n}) 
      \sum_{- \frac{q}{2} < v \le \frac{q}{2}} 
     e \left( \frac{\overline{a} v}{q} \right)
    \sigma (v, q, \gamma)
    \, 
  d \gamma  .
  %\qquad\qquad(370)
  \label{370}
\end{align}
Now we change the order of integration and summation over $v$ and use \eqref{390} to get
\begin{align} 
   \Phi''(N, d, \vec{b}) 
    &= 
  \frac{P^2}{d^4}
  \sum_{q \le P} q^{-4}
      \sum_{|\vec{n}| \le d P^{\varepsilon}} 
      e \left( \frac{n_1 b_1 + \dots + n_4 b_4}{dq} \right)
            \sum_{- \frac{q}{2} < v \le \frac{q}{2}} 
            V_q(N, d, v, \vec{b}, \vec{n})
                 \notag \\
    & \notag \\
    & \qquad \times
       \int_{  \frac{P}{2q} \le |\gamma| \le \frac{P}{q}} 
         e(- \gamma) 
  J \left( \gamma, - \frac{P}{dq} \vec{n}\right)
    \sigma (v, q, \gamma)
    \, 
  d \gamma  .
 %\qquad\qquad(380)
  \label{380}
\end{align}

%\bigskip

From \eqref{350} and \eqref{380} it follows that
\[
     \Phi''(N, d, \vec{b}) 
    \ll 
  \frac{P^2}{d^4}
  \sum_{q \le P} q^{-4}
      \sum_{|\vec{n}| \le d P^{\varepsilon}} \;
          \sum_{- \frac{q}{2} < v \le \frac{q}{2}} 
          \frac{ | V_q(N, d, v, \vec{b}, \vec{n}) |}{1 + |v|}
             \int_{  \frac{P}{2q} \le |\gamma| \le \frac{P}{q}} 
         \left|  J \left( \gamma, - \frac{P}{dq} \vec{n}\right) \right|
        \, 
  d \gamma  .
\]
Using \eqref{92} and \eqref{280} we find that 
the integral in the above formula is 
$\ll \int_{\frac{P}{2q}}^{\infty} \gamma^{-2} \, d \gamma \ll \frac{q}{P}$, hence
\begin{equation} \label{410}
     \Phi''(N, d, \vec{b}) 
    \ll 
  \frac{P}{d^4}
  \sum_{q \le P} q^{-3}
      \sum_{|\vec{n}| \le d P^{\varepsilon}} \;
          \sum_{- \frac{q}{2} < v \le \frac{q}{2}} 
          \frac{ | V_q(N, d, v, \vec{b}, \vec{n}) |}{1 + |v|} .
     %\qquad\qquad(410)
\end{equation}

%\bigskip

Next we apply the estimate for $V_q$ 
given by the inequality \eqref{420} 
from Lemma~\ref{lemma1}. 
We also notice that the sum over $v$ produces a factor $\log P$
and, having in mind that $V_q$ vanishes unless the conditions \eqref{430} hold,
we see that the summation over $\vec{n}$ produces a factor
\[
  \sum_{\substack{ |\vec{n}| \le d P^{\varepsilon} 
    \\ n_j \equiv 0 \;   ((q, d)) \\ 1 \le j \le 4
  }} 1 \ll
  \left( \frac{d P^{\varepsilon} }{ (q, d) } \right)^4 .
\]
Therefore we find
\begin{align} 
 \Phi''(N, d, \vec{b}) 
    & \ll 
    \frac{P^{1 + \varepsilon}}{d^4} \sum_{q \le P} q^{-\frac{1}{2}}
    \left( \frac{d P^{\varepsilon} }{ (q, d) } \right)^4 
    \, (q, N)^{\frac{1}{2}} \, (q, N - b_1^2 - \dots - b_4^2)^{\frac{1}{2}}
    (q, d^2)^2
         \notag \\
     & \notag \\
     & \ll    
      P^{1 + \varepsilon}   
     \sum_{q \le P} q^{-\frac{1}{2}}  \, (q, N)^{\frac{1}{2}} \, (q, N - b_1^2 - \dots - b_4^2)^{\frac{1}{2}} 
      \notag \\
     & \notag \\
     & \ll    
       P^{1 + \varepsilon}   
     \sum_{q \le P} q^{-\frac{1}{2}} \, \left(q, N (N - b_1^2 - \dots - b_4^2) \right)  .
     %\qquad\qquad(439)
      \label{439}
\end{align}

%\bigskip

For any positive integer $M$ we have
\begin{equation} \label{450}
  \sum_{q\le P} q^{-\frac{1}{2}} (q, M) \ll P^{\frac{1}{2}} \tau(M)
     %\qquad\qquad(450)
\end{equation}
(we leave the easy proof to the reader).
From the conditions  
\eqref{182} and \eqref{148} imposed on $d$ and $b_j$ we find that
$N - b_1^2 - \dots - b_4^2 \in \mathbb N$. 
Hence using \eqref{439} and \eqref{450} we find
\begin{equation} \label{460}
   \Phi''(N, d, \vec{b}) \ll P^{\frac{3}{2} + \varepsilon} .
     %\qquad\qquad(460)
\end{equation}

%\bigskip

Consider now $\Phi'(N, d, \vec{b})$. 
We remind that the expression for it is similar to the expression in the right-hand side of 
\eqref{300}, but the integration is taken over the interval 
$\left[- \frac{P}{2q}, \frac{P}{2q} \right]$.
We have
\begin{equation} \label{470}
   \Phi'(N, d, \vec{b}) = \Phi_0(N, d, \vec{b}) + \Phi^*(N, d, \vec{b}) ,
     %\qquad\qquad(470)
\end{equation}
where $\Phi_0$ denotes the contribution of the terms with $\vec{n} = \vec{0}$, that is
\begin{align}
  \Phi_0(N, d, \vec{b}) 
   & = 
   \frac{P^2}{d^4}
  \sum_{q \le P} q^{-4} \sum_{ a \; (q)^* }
  e \left( \frac{a(b_1^2 + \dots + b_4^2 -N)}{q} \right)
      \notag \\
     & \notag \\
    & \qquad \qquad \times
  G(q, ad^2, 2ad \vec{b} ) \int_{|\gamma| \le \frac{P}{2q}} e(- \gamma) 
  J \left( \gamma, \vec{0}\right) \, d \gamma  .
   %\qquad\qquad(480)
   \label{480} 
\end{align}
Respectively, $\Phi^*$ is the contribution coming from the other terms:
\begin{align}
    \Phi^*(N, d, \vec{b}) 
   & =
     \frac{P^2}{d^4}
  \sum_{q \le P} q^{-4} \sum_{a \; (q)^* }
  e \left( \frac{a(b_1^2 + \dots + b_4^2 -N)}{q} \right)
    \sum_{1 \le |\vec{n}| \le d P^{\varepsilon}} 
      e \left( \frac{n_1 b_1 + \dots + n_4 b_4}{dq} \right)
    \notag \\
    & \notag \\
    & \qquad \qquad \times
  G(q, ad^2, 2ad \vec{b} + \vec{n}) \int_{|\gamma| \le \frac{P}{2q}} e(- \gamma) 
  J \left( \gamma, - \frac{P}{dq} \vec{n}\right) \, d \gamma  .
   %\qquad\qquad(490)
    \label{490}
\end{align}

%\bigskip

Consider first $\Phi^*$. Using \eqref{390} and \eqref{490} we find
\[
   \Phi^*(N, d, \vec{b}) \ll
   \frac{P^2}{d^4}
  \sum_{q \le P} q^{-4}   
    \sum_{1 \le |\vec{n}| \le d P^{\varepsilon}  }
      | V_q(N, d, 0, \vec{b}, \vec{n}) |
   \int_{|\gamma| \le \frac{P}{2q}}
  \left|
    J \left( \gamma, - \frac{P}{dq} \vec{n}\right)
    \right|
     \, d \gamma  .
\]
We apply \eqref{93.5} to get
\[
 \int_{|\gamma| \le \frac{P}{2q}}
  \left|
    J \left( \gamma, - \frac{P}{dq} \vec{n}\right)
    \right|
     \, d \gamma 
     \ll \left( \frac{P}{qd} \, |\vec{n}| \right)^{-1 + \varepsilon} ,
\]
hence
\[
   \Phi^*(N, d, \vec{b}) \ll
   \frac{P^{1 + \varepsilon}}{d^3}
  \sum_{q \le P} q^{-3}   
     \sum_{1 \le |\vec{n}| \le d P^{\varepsilon}  }
     \frac{ | V_q(N, d, 0, \vec{b}, \vec{n}) | }{|\vec{n}|}   .
\]
Now we apply Lemma~\ref{lemma1} and find
\begin{equation} \label{500}
   \Phi^*(N, d, \vec{b}) \ll
   \frac{P^{1 + \varepsilon}}{d^3}
  \sum_{q \le P} 
  \frac{ (q, N)^{\frac{1}{2}} \, (q, N-b_1^2 - \dots - b_4^2)^{\frac{1}{2}} \, (q, d^2)^2 }{q^{\frac{1}{2}}}
     \sum_{\substack{1 \le |\vec{n}| \le d P^{\varepsilon}  \\ n_j \equiv 0 \; ((q, d)) \\ 1 \le j \le 4 }}
     \frac{1}{|\vec{n}|}   .
     %\qquad\qquad(500)
\end{equation}
It is clear that the sum over $\vec{n}$ in the expression above is
\[
  \ll \sum_{\substack{1 \le n_4 \le d P^{\varepsilon} \\ n_4 \equiv 0 ((q, d))}} \frac{1}{n_4} \;
  \left( 
  \sum_{\substack{|n|  \le n_4  \\ n \equiv 0 ((q, d)) }} 1 
  \right)^3
  \ll \sum_{1 \le h \le \frac{d P^{\varepsilon}}{(q, d )}}
  \frac{h^2}{(q, d)} \ll \frac{d^3 P^{\varepsilon}}{(q, d)^4} ,
\]
which, together with \eqref{500}, gives
\begin{align}
  \Phi^*(N, d, \vec{b}) 
  & 
  \ll
  P^{1+ \varepsilon} 
  \sum_{q \le P} \frac{(q, N)^\frac{1}{2} (q, N- b_1^2 - \dots - b_4^2)^{\frac{1}{2}}}{q^{\frac{1}{2}}} 
  \notag \\
  & \notag \\
  &
  \ll 
   P^{1+ \varepsilon} 
  \sum_{q \le P} \frac{ \left(q, N(N- b_1^2 - \dots - b_4^2) \right)}{q^{\frac{1}{2}}} .
  \notag
\end{align}
The last expression 
coincides with the expression in \eqref{439}, so we get
\begin{equation} \label{510}
   \Phi^*(N, d, \vec{b}) \ll
   P^{\frac{3}{2} + \varepsilon} .
     %\qquad\qquad(510)
\end{equation}

%\bigskip

Consider now the quantity $\Phi_0(N, d, \vec{b})$, defined by \eqref{480}. 
We use \eqref{280} and \eqref{390} to write it in the form
\[
  \Phi_0(N, d, \vec{b}) = \frac{P^2}{d^4}
  \sum_{q \le P} \frac{V_q(N, d, 0, \vec{b}, \vec{0})}{q^4} 
  \int_{|\gamma| \le \frac{P}{2q}} e(- \gamma) J(\gamma, \vec{0}) \, d \gamma .
\]
Using the estimate \eqref{92} we find that the integral in the above formula is equal to $\varkappa + O \left( \frac{q}{P} \right)$, where $\varkappa$ is defined by \eqref{90}.
We combine this with the estimate for $V_q$, given in Lemma~\ref{lemma1}, and working as above we get
\[
  \Phi_0(N, d, \vec{b}) = 
  \varkappa \frac{P^2}{d^4} \sum_{q \le P}
  \frac{V_q(N, d, 0, \vec{b}, \vec{0})}{q^4} + O \left( P^{\frac{3}{2} + \varepsilon} \right).
\]

%\bigskip

Now we extend the summation over $q$ to infinity. Using again Lemma~\ref{lemma1} and the estimate
\[
 \sum_{q > P} \frac{(q, N)^{\frac{1}{2}} (q, N-b_1^2- \dots - b_4^2)^{\frac{1}{2}} }{q^{\frac{3}{2}- \varepsilon}}
  \ll 
     \sum_{q > P} \frac{\left(q, N(N-b_1^2- \dots - b_4^2 ) \right) }{q^{\frac{3}{2}- \varepsilon}}
  \ll 
   P^{-\frac{1}{2} + \varepsilon} 
\]
(we leave the details to the reader), we find
\begin{equation} \label{520}
   \Phi_0(N, d, \vec{b}) = 
    \varkappa \frac{P^2}{d^4}  \, \sigma(N, d, \vec{b}) + 
   O \left( P^{\frac{3}{2} + \varepsilon} \right) ,
     %\qquad\qquad(520)
\end{equation}
where
\begin{equation} \label{530}
  \sigma(N, d, \vec{b}) = \sum_{q=1}^{\infty} A_q(N, d, \vec{b}) , \qquad
    A_q(N, d, \vec{b}) =  \frac{V_q(N, d, 0, \vec{b}, \vec{0})}{q^4} .
     %\qquad\qquad(530)
\end{equation}

%\bigskip

From \eqref{295}, \eqref{340}, \eqref{460}, \eqref{470}, \eqref{510} and \eqref{520}
we obtain
\begin{equation} \label{532}
    \Phi(N, d, \vec{b}) = 
      \varkappa \frac{P^2}{d^4}  \, \sigma(N, d, \vec{b}) + 
   O \left( P^{\frac{3}{2} + \varepsilon} \right) .
     %\qquad\qquad(532)
\end{equation}

%\bigskip

Now we use \eqref{170} and \eqref{532} to get
\begin{equation} \label{534}
    F(N, d) = 
      \varkappa \frac{P^2}{d^4}  \, \mathcal H (N, d) + 
   O \left( \mathcal L(N, d) \, P^{\frac{3}{2} + \varepsilon} \right) ,
     %\qquad\qquad(534)
\end{equation}
where
\begin{equation} \label{536}
  \mathcal H (N, d) = \sum_{ \vec{b} \in \mathbb Z^4 \; : \; \eqref{182} }
  \sigma(N, d, \vec{b}) .
     %\qquad\qquad(536)
\end{equation}

%\bigskip

It remains to prove that 
\begin{equation} \label{538}
  \mathcal H (N, d) = d^4 a(N) \Psi (N, d) ,
     %\qquad\qquad(538)
\end{equation}
where $a(N)$ and $\Psi(N, d)$ are defined respectively by 
\eqref{100} and \eqref{120}.
If we establish this identity and use \eqref{85}, \eqref{80} and \eqref{534} we obtain
\eqref{149} 
and finish the proof of Lemma~\ref{theorem1.5}.

%\bigskip

To prove  \eqref{538} we find an explicit formula for $\sigma(N, d, \vec{b})$.
We have already established that the series in \eqref{530} is absolutely convergent.
Further, the function $A_q(N, d, \vec{a})$ is multiplicative
with respect to $q$. Indeed, from \eqref{540} we find that
if $(q', q'')=1$ then
\[
  V_{q' q''} (N, d, 0, \vec{b}, \vec{0}) =
  V_{q'} (N,  q'' d, 0, \vec{b}, \vec{0}) \, V_{q''} (N, q' d, 0, \vec{b}, \vec{0}) .
\]
However it is easy to see that
$V_{q'} (N,  q'' d, 0, \vec{b}, \vec{0}) = V_{q'} (N,  d, 0, \vec{b}, \vec{0})$
and $V_{q''} (N, q' d, 0, \vec{b}, \vec{0}) = V_{q''} (N, d, 0, \vec{b}, \vec{0})$
and it remains to apply \eqref{530}.

%\bigskip

Hence we have
\begin{equation} \label{700}
  \sigma(N, d, \vec{b}) = 
  \prod_{p} \chi_p(N, d, \vec{b}) ,
     %\qquad\qquad(700)
\end{equation}
where
\begin{equation} \label{710}
    \chi_p(N, d, \vec{b}) = 1 + \sum_{s=1}^{\infty}  A_{p^s}(N, d, \vec{b}) .
     %\qquad\qquad(710)
\end{equation}

%\bigskip

We shall now compute the quantities $\chi_p(N, d, \vec{b})$.

%\bigskip

Consider first the case $p \nmid 2d$. A straightforward calculation, based on 
\eqref{253}, \eqref{252}, \eqref{260}, \eqref{545.6}, \eqref{390} and \eqref{530},
shows that
\[
  A_{p^s}(N, d, \vec{b}) = p^{-2s} c_{p^s} (N) ,
\]
Using \eqref{545.7} and \eqref{110} we find that 
\begin{equation} \label{720}
  A_{p^s}(N, d, \vec{b}) = 
  \begin{cases}
     0 \quad & \text{for} \quad s \ge \xi_p(N) + 2 , \\
     - \frac{1}{p^{\xi_p(N) + 2}}  & \text{for} \quad s = \xi_p(N) + 1 , \\
     \frac{1}{p^s} - \frac{1}{p^{s+1}} & \text{for} \quad s \le \xi_p(N) . 
  \end{cases}
     %\qquad\qquad(720)
\end{equation}
(The third case in \eqref{720} is applicable only if $\xi_p(N) \ge 1$).

%\bigskip

From \eqref{710} and \eqref{720} we find
\begin{equation} \label{730}
  \chi_{p}(N, d, \vec{b}) = 
    \left( 1 + \frac{1}{p} \right)
    \left( 1 - \frac{1}{p^{\xi_p(N) + 1}} \right)
       \qquad \text{for} \qquad  p \nmid 2 d .
     %\qquad\qquad(730)
\end{equation}

%\bigskip

Suppose now that $p \mid d$. From \eqref{182} it follows that $p \mid N - b_1^2 - \dots - b_4^2$
and using 
\eqref{251}, \eqref{260}, \eqref{390} and \eqref{530}, 
we easily find that
\begin{equation} \label{731}
  A_p(N, d, \vec{b}) = p-1 \qquad \text{for} \qquad p \mid d .
     %\qquad\qquad(731)
\end{equation}

%\bigskip

Suppose now that $s \ge 2$. In this case we have $(p^s, a d^2)=p^2$. However 
$p \nmid b_j$ for any $j$ because of the condition $d \mid b_1b_2b_3b_4 + 1$ imposed in \eqref{182}. 
This implies that $p^2 \nmid 2adb_j$ and using \eqref{251} we find
$G(p^s, a d^2, 2adb_j) = 0$. Therefore from \eqref{390}
and \eqref{530} we find
\begin{equation} \label{732}
  A_{p^s}(N, d, \vec{b}) = 0  \qquad \text{for} \qquad s \ge 2 , \qquad p \mid d .
     %\qquad\qquad(732)
\end{equation}
From \eqref{710}, \eqref{731} and \eqref{732}
we obtain
\begin{equation} \label{740}
  \chi_{p}(N, d, \vec{b}) 
   =  p  \qquad  \text{for} \qquad p \mid d .
        %\qquad\qquad(740)
\end{equation}

%\bigskip

It remains to consider the case $p=2$. Using \eqref{250}, \eqref{254}, \eqref{255}, \eqref{545.6}, \eqref{545.7},
\eqref{390},  \eqref{530} and our assumption $2 \nmid N$ we easily get
\[
    A_{2^s} (N, d, \vec{b}) = 0 \qquad \text{for} \qquad  s \ge 1
\]
(we leave the verification to the reader). 
This formula and \eqref{710} imply
\begin{equation} \label{750}
  \chi_{2}(N, d, \vec{b}) = 1 .
        %\qquad\qquad(750)
\end{equation}

%\bigskip

From \eqref{700}, \eqref{730}, \eqref{740} and \eqref{750}
we get
\begin{equation} \label{752}
  \sigma(N, d, \vec{b}) = 
      d \,  \prod_{ p \nmid 2d} \left( 1 + \frac{1}{p} \right)
      \left( 1 - \frac{1}{p^{\xi_p(N)+1}} \right) 
        %\qquad\qquad(752)
\end{equation}
and bearing in mind the definitions \eqref{100} and \eqref{130} we obtain
\begin{equation} \label{800}
  \sigma(N, d, \vec{b}) = 
       d \, a(N) \, \alpha(N, d) .
     %\qquad\qquad(800)
\end{equation}

%\bigskip

From \eqref{120}, \eqref{536} and \eqref{800} we find that the quantity
$\mathcal H(N, d)$ satisfies \eqref{538}
and the proof of Lemma~\ref{theorem1.5} is complete.

\hfill
$\square$

\subsection{The main term}

\indent

To apply the sieve method and prove the theorem we have to study the properties of the 
main term $M(N, d)$ defined by \eqref{85}. We already mentioned that the constant 
$\varkappa$ satisfies \eqref{93}. Further, from \eqref{100} we easily find
\begin{equation} \label{790}
  1 \ll a(N) \ll \log \log N .
     %\qquad\qquad(790)
\end{equation}

%\bigskip

More care is needed about the quantity $\Psi(N, d)$ defined by \eqref{120}.
We have the following 
\begin{lemma} \label{lemma791}
The function $\Psi(N, d)$ is multiplicative with respect to $d$.
We also have 
\begin{equation} \label{791}
    \Psi(N, p) < 0.9 \qquad \text{for} \qquad p > 2
     %\qquad\qquad(791)
\end{equation}
and
\begin{equation} \label{792}
   0 < \Psi(N, p)  \qquad \text{for} \qquad p > 1000 .
     %\qquad\qquad(792)
\end{equation}
Finally, for all $z_1, z_2$ with $ 2 < z_1 < z_2$ we have
\begin{equation} \label{792.5}
   \prod_{z_1 \le p < z_2} \left( 1 - \Psi(N, p) \right)^{-1} \le
     \frac{\log z_2}{ \log z_1}
        \left( 1 + \frac{L}{\log z_1} \right) ,
     %\qquad\qquad(792.5)
\end{equation}
where $L > 0$ is an absolute constant.
\end{lemma}

\paragraph{Proof.}
Obviously the function $\alpha(N, d)$, defined by \eqref{130} is multiplicative 
with respect to $d$ and it is easy to see that the same property possesses 
$\mathcal L(N, d)$, which, by definition, is the number of solutions of the system \eqref{140}.
This proves the multiplicativity of $\Psi(N, d)$.

%\bigskip

We shall now study $\Psi(N, p)$ for $p>2$.

%\bigskip

It is easy to verify that for any prime $p > 2$ we have
\begin{equation} \label{793}
    \frac{3}{4} \le \alpha(N, p) \le \frac{9}{8} .
     %\qquad\qquad(793)
\end{equation}

%\bigskip

Consider $\mathcal L = \mathcal L(N, p)$.
We shall prove that for $p > 2$ we have
\begin{equation} \label{793.5}
    \mathcal L \le 4 (p-1)^2
     %\qquad\qquad(793.5)
\end{equation}
and
\begin{equation} \label{794}
   \left| \mathcal L - p^2 \right| \le 30 p^{\frac{3}{2}} .
     %\qquad\qquad(794)
\end{equation}

%\bigskip

Suppose that the integers $b_1, \dots, b_4$ satisfy 
\[
  b_1^2 + \dots + b_4^2 \equiv N \pmod{p} , \qquad b_1b_2b_3b_4 + 1 \equiv 0 \pmod{p} .
\]
From the second of these congruences we conclude that
$b_4 \equiv - \overline{b_1 b_2 b_3} \pmod{p}$, hence
$\mathcal L $ is equal to the number of triples $b_1, b_2, b_3 \in \{ 1, 2, \dots , p-1 \}$
satisfying 
\[
   b_1^2 + b_2^2 + b_3^2 + \overline{b_1 b_2 b_3}^{\, 2} \equiv N \pmod{p} ,
\]
or equivaliently
\begin{equation} \label{900}
   b_1^2 b_2^2 b_3^2 \left( b_1^2 + b_2^2 + b_3^2  - N \right) + 1  \equiv 0  \pmod{p} .
     %\qquad\qquad(900)
\end{equation}
For fixed $b_1, b_2$ there are at most $4$ admissible values of $b_3$, hence 
\eqref{793.5} is correct.

%\bigskip

Using the definition of $\Psi(N, p)$ given by \eqref{120} as well as \eqref{793}, \eqref{793.5}
we establish \eqref{791}.

%\bigskip

To establish \eqref{792} we first prove \eqref{794} in the following elementary way.
For any integer $a$ the number of solutions of 
$x^2 \equiv a \pmod{p}$ is equal to $1 + \left( \frac{a}{p} \right)$, so we may write
\[
 \mathcal L = \sum_{\substack{a_1, a_2, a_3 \; (p)^* \\ \eqref{902} }}
   \left( 1 + \left( \frac{a_1}{p} \right) \right)
   \left( 1 + \left( \frac{a_2}{p} \right) \right)
   \left( 1 + \left( \frac{a_3}{p} \right) \right) ,
\]
where the summation is taken over variables $a_1, a_2, a_3$ satisfying
\begin{equation} \label{902}
   a_1 a_2 a_3 \left( a_1 + a_2 + a_3  - N \right) + 1  \equiv 0  \pmod{p} .
     %\qquad\qquad(902)
\end{equation}
It is clear that
\begin{equation} \label{904}
  \mathcal L = \mathcal L_1 + 3 \mathcal L_2 + 3 \mathcal L_3 + \mathcal L_4 ,
     %\qquad\qquad(904)
\end{equation}
where $\mathcal L_1$ is the number of solutions of \eqref{902},
\begin{equation} \label{906}
  \mathcal L_2
    = 
   \sum_{\substack{a_1, a_2, a_3 \; (p)^* \\ \eqref{902} }}
   \left( \frac{a_1}{p} \right)   ,
   \quad
    \mathcal L_3
   = \sum_{\substack{a_1, a_2, a_3 \; (p)^* \\ \eqref{902} }}
   \left( \frac{a_1 a_2}{p} \right)  ,
   \quad
    \mathcal L_4
    = 
   \sum_{\substack{a_1, a_2, a_3 \; (p)^* \\ \eqref{902} }}
   \left( \frac{a_1 a_2 a_3}{p} \right)  .
   %\qquad\qquad(906)
\end{equation}

%\bigskip

Consider $\mathcal L_4$. 
We use the identity
\begin{equation} \label{906.5}
        \sum_{h \; (p)} e \left( \frac{mh}{p} \right) 
        = 
         \begin{cases}
         p \qquad & \text{if} \qquad p \mid m , \\
         0  & \text{otherwise} 
         \end{cases}
   %\qquad\qquad(906.5)
\end{equation}
and find
\[
  \mathcal L_4 = \frac{1}{p} 
  \sum_{h, a_1, a_2, a_3 \; (p)^*}
  \left( \frac{a_1 a_2 a_3}{p} \right) 
  e \left( \frac{h \left( a_1 a_2 a_3 (a_1 + a_2 + a_3 - N ) + 1 \right)}{p} \right) .
\]
For fixed $a_2, a_3, h$ we change the variable $a_1$ to $a_4$, where
$a_1 a_2 a_3 \equiv a_4 \pmod{p}$. We find
\[
  \mathcal L_4 = \frac{1}{p} 
  \sum_{h, a_2, a_3, a_4 \; (p)^*}
  \left( \frac{a_4}{p} \right) 
  e \left( \frac{h \left( a_4 (\overline{a_2 a_3} a_4 + a_2 + a_3 - N ) + 1 \right)}{p} \right) .
\]
Next we change the variable $h$ to $t$, where $h \equiv a_2 t \pmod{p}$, and use \eqref{250} to get 
\begin{align}
  \mathcal L_4 
   & = 
   \frac{1}{p} 
  \sum_{t, a_2, a_3, a_4 \; (p)^*}
  \left( \frac{a_4}{p} \right) 
    e \left( \frac{ t \left( \overline{a_3} a_4^2 + a_4 a_2^2 + (a_3 a_4 - a_4 N + 1 )a_2 \right)}{p} \right) 
    \notag \\
    & \notag \\
   & =
       \frac{1}{p} 
  \sum_{t, a_3, a_4 \; (p)^*}   \left( \frac{a_4}{p} \right) 
    e \left( \frac{ t  \overline{a_3} a_4^2 }{p} \right) 
    \left( G(p, ta_4, t (a_3 a_4 - a_4 N + 1)) - 1 \right) .
    \notag
\end{align}
Clearly the contribution of the term $-1$ in the brackets above vanishes.
Therefore, using \eqref{252} we find
\[
  \mathcal L_4 =
   \frac{G(p, 1)}{p}  \sum_{t, a_3, a_4 \; (p)^*} 
    \left( \frac{t}{p} \right) 
     e \left( \frac{ t \left(  \overline{a_3} a_4^2 - \overline{4 a_4} (a_3a_4 - a_4 N + 1)^2 \right) }{p} \right) 
\]
We write the summation over $t$ inside and applying \eqref{253} and \eqref{256} we find
\begin{align}
  \mathcal L_4 
  & =
   \frac{G^2(p, 1)}{p}  \sum_{a_3, a_4 \; (p)^*} 
       \left( \frac{   \overline{a_3} a_4^2 - \overline{4 a_4} (a_3a_4 - a_4 N + 1)^2  }{p} \right) 
        \notag \\
        & \notag \\
    & =    
       (-1)^{\frac{p-1}{2}} 
       \sum_{a_4 \; (p)^*} \left( \frac{a_4}{p} \right) 
       \sum_{a_3 \; (p)}    
        \left( \frac{   4a_3 a_4^3 - a_3^2 (a_3a_4 - a_4 N + 1)^2  }{p} \right) .
       \notag
\end{align}
We estimate the character sum over $a_3$ using \eqref{27} and find that its modulus does not exceed 
$3 \sqrt{p}$. Hence we obtain
\begin{equation} \label{910}
  | \mathcal L_4 | \le 3 p^{\frac{3}{2}} .
     %\qquad\qquad(910)
\end{equation}

%\bigskip

Consider now $\mathcal L_3$, From \eqref{250}, \eqref{902}, \eqref{906} and \eqref{906.5} we find
\begin{align}
   \mathcal L_3 
    & = 
   \frac{1}{p} \sum_{h, a_1, a_2, a_3 \; (p)^*} 
     \left( \frac{a_1 a_2}{p} \right) 
       e \left( \frac{h \left( a_1 a_2 a_3 (a_1 + a_2 + a_3 - N) + 1 \right)}{p} \right)
       \notag \\
     & =
      \frac{1}{p} \sum_{h, a_1, a_2 \; (p)^*} 
     \left( \frac{a_1 a_2}{p} \right) 
     e \left( \frac{h}{p} \right)
     \left( G(p, ha_1 a_2, h a_1 a_2 (a_1 + a_2 - N)) - 1 \right) .
        \notag
\end{align}
Clearly the contribution of the term $-1$ in the brackets above vanishes.
Now we apply \eqref{252} to find
\[
    \mathcal L_3 =
    \frac{G(p, 1)}{p}  \sum_{h, a_1, a_2 \; (p)^*} 
    \left( \frac{h}{p} \right) 
     e \left( \frac{h \left( 1 - \overline{4} a_1 a_2 (a_1 + a_2  - N)^2 \right)}{p} \right) .
\]
We insert the summation over $h$ inside and use \eqref{256} to get 
\[
    \mathcal L_3 = 
     \frac{G^2(p, 1)}{p}  \sum_{a_1, a_2 \; (p)^*} 
        \left( \frac{ 1 - \overline{4} a_1 a_2 (a_1 + a_2  - N)^2 }{p} \right) .
\]
Applying \eqref{27} for the sum over $a_2$ and having also in mind \eqref{253} we obtain
\begin{equation} \label{920}
  |\mathcal L_3| \le 3 p^{\frac{3}{2}} .
     %\qquad\qquad(920)
\end{equation}

%\bigskip

In the same manner we consider $\mathcal L_2$ and find
\begin{equation} \label{921}
   |\mathcal L_2| \le 3 p^{\frac{3}{2}} .
     %\qquad\qquad(921)
\end{equation}

%\bigskip

It remains to study $\mathcal L_1$. We use \eqref{906.5} and find
\begin{equation} \label{922}
   \mathcal L_1 = \frac{1}{p} \sum_{a_1, a_2, a_3 \; (p)^*} \;
   \sum_{h \; (p)} 
   e \left( \frac{ h \left( a_1 a_2 a_3 (a_1 + a_2 + a_3 - N) + 1 \right) }{p} \right) 
   = \frac{(p-1)^3}{p} + \Delta ,
     %\qquad\qquad(922)
\end{equation}
where
\[
  \Delta 
   = 
   \frac{1}{p} \sum_{h, a_1, a_2, a_3 \; (p)^*} 
      e \left( \frac{ h \left( a_1 a_2 a_3 (a_1 + a_2 + a_3 - N) + 1 \right) }{p} \right) 
\]
Now we apply \eqref{250}, \eqref{252} and \eqref{256} to get
\begin{align}
  \Delta 
   & = 
   \frac{1}{p} \sum_{h, a_1, a_2 \; (p)^*} 
      e \left( \frac{ h }{p} \right) 
      \left( G (p, h a_1 a_2, h a_1 a_2 (a_1 + a_2 - N )) - 1 \right)
      \notag \\
      & \notag \\
   & =
         \frac{G(p, 1)}{p} \sum_{h, a_1, a_2 \; (p)^*} 
         \left( \frac{h a_1 a_2}{p} \right)
         e \left( \frac{h \left( 1 - \overline{4} a_1 a_2 (a_1 + a_2 - N)^2 \right) }{p} \right)
         + \frac{(p-1)^2}{p}
          \notag \\
      & \notag \\
   & =
    \frac{G^2(p, 1)}{p} \sum_{a_1, a_2 \; (p)^*} 
      \left( \frac{a_1 a_2 \left( 4 -  a_1 a_2 (a_1 + a_2 - N)^2  \right) }{p} \right) 
    + \frac{(p-1)^2}{p} .
      \notag
\end{align}
We estimate the sum over $a_2$ using \eqref{27} and having in mind \eqref{253}
we get
\begin{equation} \label{960}
   |\Delta| \le 4 p^{\frac{3}{2}}  .
     %\qquad\qquad(960)
\end{equation}
From \eqref{904}, \eqref{910} -- \eqref{960} we obtain \eqref{794}.

%\bigskip

The inequality \eqref{792} for $\Psi(N, p)$ follows from
\eqref{120}, \eqref{793} and \eqref{794}.

%\bigskip

It remains to prove \eqref{792.5}. 
From \eqref{120}, \eqref{130} and \eqref{794} we see that 
\begin{equation} \label{963}
  \Psi(N, p) = \frac{p^2 + O \left( p^{\frac{3}{2}} \right) }{p^3 \left( 1 + \frac{1}{p} \right) 
  \left( 1 - \frac{1}{p^{1 + \xi_p(N)}} \right) } 
  = \frac{1}{p} + O \left( \frac{1}{p^{\frac{3}{2}}} \right) ,
     %\qquad\qquad(963)
\end{equation}
with an absolute constants in the $O$-terms.
We apply Mertens's prime number theorem (see \cite[Ch. 22]{Hardy-Wright})
and after some simple calculations, which we leave to the reader, we  establish \eqref{792.5}. 

\hfill
$\square$

\subsection{End of the proof of Theorem~\ref{theorem1}}

\indent

Here we use the terminology and results from \cite[Ch. 12]{FriedIw}.

%\bigskip

We write the quantity $P(z)$ given by \eqref{21} in the form
\begin{equation} \label{1100}
  P(z) = C_0 \, P^*(z) ,
     %\qquad\qquad(1100)
\end{equation}
where
\begin{equation} \label{1110}
  C_0 = \prod_{2 < p < 1000} p , \qquad
  P^*(z) = \prod_{1000 < p < z} p .
     %\qquad\qquad(1110)
\end{equation}
Suppose that 
\begin{equation} \label{1120}
  D = N^{\delta} , \qquad 0 < \delta < \frac{1}{12} 
     %\qquad\qquad(1120)
\end{equation}
and let $\lambda(d)$ be the lower bound Rosser weights of level $D$, 
hence
\begin{equation} \label{1125}
  |\lambda(d)| \le 1 ; \qquad 
  \lambda(d) = 0 \quad \text{for} \quad d > D \quad \text{or} \quad \mu^2(d) = 0 .
     %\qquad\qquad(1125)
\end{equation}
Then for the sum $\Gamma$, defined by \eqref{50}, we have
\begin{align}
  \Gamma 
    & = 
  \sum_{ x_1^2 + x_2^2 + x_3^2 + x_4^2 = N  }
            \omega(x_1) \dots \omega(x_4)
            \sum_{\delta \mid  (x_1 x_2 x_3 x_4 + 1 , C_0)  } \mu(\delta)
                 \sum_{ t \mid  (x_1 x_2 x_3 x_4 + 1 , P^*(z)) } \mu (t)  
            \notag \\
       & \notag \\
       & \ge
       \sum_{ x_1^2 + x_2^2 + x_3^2 + x_4^2 = N  }
              \omega(x_1) \dots \omega(x_4)
                \sum_{\delta \mid  (x_1 x_2 x_3 x_4 + 1 , C_0)  } \mu(\delta)
                    \sum_{ t \mid  (x_1 x_2 x_3 x_4 + 1 , P^*(z)) } \lambda (t)   .
            \notag
\end{align}
Now we change the order of summation and find
\begin{equation} \label{1130}
  \Gamma \ge
  \sum_{d \mid P(z)} \theta (d) \, F(N, d) , 
     %\qquad\qquad(1130)
\end{equation}
where $F(N, d)$ is defined by \eqref{70} and
\begin{equation} \label{1140}
  \theta (d) = \sum_{\substack{ \delta \mid C_0 \\ t \mid P^*(z) \\ \delta t = d  } }\mu(\delta) \lambda(t) .
     %\qquad\qquad(1140)
\end{equation}
We apply \eqref{80} and \eqref{1130} and find that
\begin{equation} \label{1150}
  \Gamma \ge \Gamma_1 + R ,
     %\qquad\qquad(1150)
\end{equation}
where
\begin{equation} \label{1160}
  \Gamma_1 = \sum_{d \mid P(z)} \theta(d) \, M(N, d) , \qquad R = \sum_{d \mid P(z)} \theta(d) \, R(N, d) .
     %\qquad\qquad(1160)
\end{equation}

%\bigskip

From \eqref{1125} and \eqref{1140} we see that
$\theta(d) \ll 1$ and also that
$\theta(d)$ is supported on the set of squarefree odd integers $d \le  C_0 D $.
Therefore using Lemma~\ref{theorem1.5} we get
\[
   R \ll  \sum_{\substack{d \le C_0 D \\ 2 \nmid d}} \mu^2(d) \, |R(N, d)| 
   \ll 
   N^{\frac{3}{4} + \varepsilon} 
   \sum_{\substack{d \le C_0 D \\ 2 \nmid d}} \mu^2(d) \, \mathcal L(N, d) .
\]   
Having in mind \eqref{794} we see that for any squarefree odd $d$ we have
\[
  \mathcal L(N, d) \le d^2 \prod_{p \mid d} \left( 1 + 30 p^{-\frac{1}{2}} \right) \ll d^2 \tau(d) .
\]
Hence using \eqref{1120} we get
\begin{equation} \label{1170}   
   R \ll
     D^3 N^{\frac{3}{4} + \varepsilon} \ll \frac{N}{\log^2 N}  .
     %\qquad\qquad(1170)
\end{equation}

%\bigskip

Consider now the sum $\Gamma_1$. Using \eqref{85} we write it as
\begin{equation} \label{1180}
  \Gamma_1 = \varkappa \, N \, a(N) \, \Gamma_2 ,
     %\qquad\qquad(1180)
\end{equation}
where
\[
  \Gamma_2 = \sum_{d \mid P(z)} \theta(d) \, \Psi(N, d) .
\]
Using the muptiplicativity with respect to $d$ of $\Psi(N, d)$ and having in mind
\eqref{1140} we find
\begin{equation} \label{1190}
  \Gamma_2 = \sum_{\substack{ \delta \mid C_0 \\ t \mid P^*(z)}}
   \mu(\delta) \lambda(t) \Psi(N, \delta t) = 
   \Gamma_3 \, \Gamma_4 ,
     %\qquad\qquad(1190)
\end{equation}
where
\[
  \Gamma_3 
     = 
    \sum_{\delta \mid C_0} \mu(\delta) \, \Psi(N, \delta)  ,
     \qquad
   \Gamma_4 = 
    \sum_{t \mid P^*(z)} \lambda(t) \, \Psi(N, t)   .
\]

%\bigskip

Using \eqref{791} and \eqref{1110} we find
\begin{equation} \label{1200}
  \Gamma_3  = \prod_{2 < p < 1000} \left( 1 - \Psi(N, p)\right) \gg 1 ,
     %\qquad\qquad(1200)
\end{equation}
where the constant in Vinogradov's symbol is absolute.

%\bigskip

Consider now $\Gamma_4$. We apply the lower bound linear sieve and 
having in mind the properties of $\Psi(N, d)$ mentioned in Lemma~\ref{lemma791}
we obtain
\begin{equation} \label{1210}
  \Gamma_4 \ge \Pi (z) \left( f(s_0) + O \left( (\log D)^{- \frac{1}{3}} \right) \right)  ,
     %\qquad\qquad(1210)
\end{equation}
where
\begin{equation} \label{1220}
  \Pi (z) = \prod_{1000 < p < z} \left( 1 - \psi (N, p) \right) ,
     %\qquad\qquad(1220)
\end{equation}
\begin{equation} \label{1217}
  s_0 = \frac{\log D}{\log z} = \frac{\delta}{\eta}
     %\qquad\qquad(1217)
\end{equation}
and where $f(s)$ is the lower function of the linear sieve,
for which we know that
\begin{equation} \label{1215}
  f(s) = 2 e^{\gamma}  s^{-1}  \log(s-1)  \qquad \text{for} \qquad s \in (2, 3) 
     %\qquad\qquad(1215)
\end{equation}
($\gamma$ is the Euler constant).

%\bigskip

We choose
\begin{equation} \label{1219}
    \eta = \frac{1}{24} - 10^{-4} , \qquad \delta = \frac{1}{12} - 10^{-4} .
     %\qquad\qquad(1219)
\end{equation}
Then from \eqref{21}, \eqref{963} and \eqref{1219} it follows that 
\begin{equation} \label{1230}
  \Pi (z) \asymp (\log z)^{-1}  \asymp (\log N)^{-1} .
     %\qquad\qquad(1230)
\end{equation}
On the other hand from \eqref{1217} and \eqref{1219} we find $s_0 \in (2, 3)$ and having in mind \eqref{1215}
we find that 
\begin{equation} \label{1240}
   f(s_0) > 0 .
     %\qquad\qquad(1240)
\end{equation}
From \eqref{1120}, \eqref{1210}, \eqref{1230} and \eqref{1240} we find
\[
  \Gamma_4 \gg (\log N)^{-1}
\]
and having also in mind \eqref{790}, \eqref{1150}, \eqref{1170} -- \eqref{1200} we obtain \eqref{60}. 
It remains to notice that for the number $\eta$ given by \eqref{1219}
we have $48 < \frac{2}{\eta} < 49$ and the theorem is proved.

\hfill
$\square$

%\bigskip
%\bigskip

\vbox{
\hbox{Faculty of Mathematics and Informatics}
\hbox{Sofia University ``St. Kl. Ohridsky''}
\hbox{5 J.Bourchier, 1164 Sofia, Bulgaria}
\hbox{ }
\hbox{tlt@fmi.uni-sofia.bg}
\hbox{dtolev@fmi.uni-sofia.bg}}

\end{document}